\documentclass[12pt]{amsart}    
\usepackage{amscd}

\def\NZQ{\Bbb}               
\def\NN{{\NZQ N}}

\def\ZZ{{\NZQ Z}}

\def\frk{\frak}               

\def\mm{{\frk m}}

\def\opn#1#2{\def#1{\operatorname{#2}}} 
\opn\chara{char}
\opn\length{\ell}
\opn\pd{pd}
\opn\rk{rk}
\opn\projdim{proj\,dim}
\opn\rank{rank}
\opn\depth{depth}
\opn\grade{grade}
\opn\height{ht}
\opn\embdim{emb\,dim}
\opn\codim{codim}

\opn\Tr{Tr}
\opn\bigrank{big\,rank}
\opn\superheight{superheight}\opn\lcm{lcm}
\opn\trdeg{tr\,deg}%
\opn\reg{reg}
\opn\lreg{lreg}
\opn\div{div}
\opn\Div{Div}
\opn\WDiv{WDiv}
\opn\cl{cl}
\opn\Cl{Cl}
\opn\Spec{Spec}
\opn\Supp{Supp}
\opn\supp{supp}
\opn\Sing{Sing}
\opn\Ass{Ass}
\opn\Assh{Assh}
\opn\Min{Min}
\opn\Reg{Reg}
\opn\Ann{Ann}
\opn\Rad{Rad}
\opn\Soc{Soc}
\opn\Socle{Socle}
\opn\Ker{Ker}
\opn\Coker{Coker}
\opn\Im{Im}
\opn\Hom{Hom}
\opn\Tor{Tor}
\opn\Ext{Ext}
\opn\End{End}
\opn\Aut{Aut}
\opn\id{id}

\opn\nat{nat}
\opn\pff{pf}
\opn\Pf{Pf}
\opn\GL{GL}
\opn\SL{SL}
\opn\mod{mod}
\opn\ord{ord}
\opn\Proj{Proj}
\opn\aff{aff}
\opn\con{conv}
\opn\relint{relint}
\opn\st{st}
\opn\lk{lk}
\opn\cn{cn}
\opn\core{core}
\opn\vol{vol}
\opn\link{link}
\opn\star{star}
\opn\gr{gr}

\def\pot#1#2{#1[\kern-0.28ex[#2]\kern-0.28ex]}

\opn\dirlim{\underrightarrow{\lim}}
\opn\inivlim{\underleftarrow{\lim}}
\let\iso=\cong
\let\Union=\bigcup
\let\Sect=\bigcap
\let\Dirsum=\bigoplus

\let\mcone= * 
\let\To=\longrightarrow
\def\Implies{\ifmmode\Longrightarrow \else
     \unskip${}\Longrightarrow{}$\ignorespaces\fi}
\def\implies{\ifmmode\Rightarrow \else
     \unskip${}\Rightarrow{}$\ignorespaces\fi}
\def\iff{\ifmmode\Longleftrightarrow \else
     \unskip${}\Longleftrightarrow{}$\ignorespaces\fi}

\let\:=\colon
\opn\H{H}
\opn\Pic{Pic}

\newtheorem{Theorem}{Theorem}

\newtheorem{Lemma}[Theorem]{Lemma}
\newtheorem{Corollary}[Theorem]{Corollary}
\newtheorem{Proposition}[Theorem]{Proposition}

\newtheorem{Remark}{Remark}

\newtheorem{Definition}{Definition}

\let\epsilon\varepsilon
\opn\inii{in}
\opn\inim{inm}
\opn\set{set}
\def\pnt{{\raise0.5mm\hbox{\large\bf.}}}

\begin{document}

\title{Local cohomologies of Isolated non $F$-rational Singularities}
\author{Yukihide Takayama}
\address{Yukihide Takayama, Department of Mathematical
Sciences, Ritsumeikan University, 
1-1-1 Nojihigashi, Kusatsu, Shiga 525-8577, Japan}
\email{takayama@se.ritsumei.ac.jp}

\def\Coh#1#2{H_{\mm}^{#1}(#2)}
\def\eCoh#1#2#3{H_{#1}^{#2}(#3)}

\newcommand{\AppTh}{Theorem~\ref{approxtheorem} }
\def\da{\downarrow}
\newcommand{\ua}{\uparrow}
\newcommand{\namedto}[1]{\buildrel\mbox{$#1$}\over\rightarrow}
\newcommand{\bdel}{\bar\partial}
\newcommand{\proj}{{\rm proj.}}

\maketitle

\newenvironment{myremark}[1]{{\bf Note:\ } \dotfill\\ \it{#1}}{\\ \dotfill
{\bf Note end.}}
\newcommand{\transdeg}[2]{{\rm trans. deg}_{#1}(#2)}
\newcommand{\mSpec}[1]{{\rm m\hbox{-}Spec}(#1)}

\newcommand{\tbf}{{{\Large To Be Filled!!}}}

\pagestyle{plain}
\maketitle
\def\gCoh#1#2#3{H_{#1}^{#2}\left(#3\right)}
\def\subsetneq{\raisebox{.6ex}{{\small $\; \underset{\ne}{\subset}\; $}}}

\begin{abstract}
In this paper, we consider positively graded isolated non $F$-rational
singularities $(R,\mm)$ with $d=\dim R$ over the field $K$
of positive characteristic.
 We give a representation of lower local
cohomologies $\Coh{i}{R}$ $(i<d)$ in terms of tight closure and
limit closure of certain type of parameters. 
As an application to isolated
singularities, we show a relation between non-vanishing of the tight
closure of zero in the highest local cohomology $(0)^*_{\Coh{d}{R}}$
and non-vanishing of the cohomology $\Coh{d-1}{R}$.\\
MSC Primary: 13A35, 13D45, Secondary: 13A02
\end{abstract}

\section*{Introduction}
Let $K$ be a field of characteristic $\chara{K}=p >  0$.
Let $(R, \mm)$ be a finitely generated $\NN$-graded $K$-algebra with $d
= \dim R$, $R = \Dirsum_{n\geq 0}R_n$, $\mm = \Dirsum_{n>0}R_n$ 
and $R_0 = K$.
We assume that $R$ is reduced, equidimensional and $R_P$ is
$F$-rational for all primes $P(\ne\mm)$, namely $R$ is an 
$\NN$-graded isolated non $F$-rational singularity.
Since $R$ is
Cohen-Macaulay on the punctured spectrum, it is  a generalized
Cohen-Macaulay ring, i.e., $\length_R(\Coh{i}{R})<\infty$ for all $i<d$,
where $\Coh{i}{R}$ is the $i$th local cohomology module and 
$\length_R(-)$ denotes the length as a $R$-module.
In this paper, we are interested in the structure of local cohomologies
$\Coh{i}{R}$.

A distinguished property of $\NN$-graded isolated non $F$-rational
singularity is Kodaira vanishing theorem, which, in terms of local
cohomologies, can be stated that lower local cohomologies vanish for
negative degrees \cite{HuSIX,HuSm}: $[\Coh{i}{R}]_n = 0$ for all $i<d$
and all $n<0$. If $R$ is a generic mod $p$ reduction
of an isolated non-rational singularity (over the field of characteristic 
$0$), Kodaira vanishing holds. This result has been established by Huneke, 
Smith, Hara and Watanabe \cite{Hara4,HuSm,HW96}. 
This is in fact the case of {\em liftable to the second Witt vectors} 
in the sence of Deligne-Illusie \cite{DI}. 
Notice that an isolated non $F$-rational singularity
is not necessary a mod $p$ reduction of an isolated 
non-rational singularity and 
Kodaira vanishing for positive characteristic is 
in general false, eg. \cite{Rey}.
This means that vanishing/non-vanishing
 of local cohomologies for isolated non $F$-rational
singularities in general is widely open.

In this papaer, we first show that lower local cohomologies
$\Coh{i}{R}$, $i<d$, of an $\NN$-graded isolated non $F$-rational
singularity $R$ can be described as quotient modules of the tight
closures by the limit closures of suitable parameter ideals (see
Theorem~\ref{1}). This is a refinement of the results by Schenzel
\cite{Sch} and Smith \cite{KS8}.  With this representation, 
Huneke-Smiths interpretation of Kodaira vanishing in
terms of tight closure \cite{HuSm} can be recovered immediately in a
slightly general form.

A Noetherian local ring $(R,\mm)$ is an isolated non $F$-rational
singularity if and only if $R$ is generalized Cohen-Macaulay and 
the tight closure of zero in the highest local cohomology
$(0)^*_{\Coh{d}{R}}$, $d=\dim R$, has finite length (cf. \cite{GN02}).
This suggests that vanishing/non-vanishing of lower local 
cohomologies may be controlled by $(0)^*_{\Coh{d}{R}}$ to some 
extent. 
In fact, the above mensioned 
Kodaira vanishing theorem 
by Huneke-Smith-Hara-Watanabe has been proved by investigating 
the structure of $(0)^*_{\Coh{d}{R}}$.
In this paper, we will consider isolated 
singularities and show that $[(0)^*_{\Coh{d}{R}}]_n\ne 0$
for some $n\in\ZZ$ implies non-vanishing of 
$\Coh{d-1}{R}$ for certain degree, under some subtle conditions
(see Theorem~\ref{main}).

In section~1 we will summerize the known results on isolated 
non $F$-rational singularities that are necessary in this paper.
In section~2, after preparing some definitions and facts 
on standard sequence, limit closure and germ closure, we 
give a representation of lower local cohomologies 
in terms of tight closure and limit closure. We also give 
some direct consequences from this representation.
$\NN$-graded isolated singularities are considered in section~3.
After giving a graded version of Goto-Nakamura's representation
of tight closure of zero $(0)^*_{\Coh{d}{R}}$, we 
investigate non-vanishing of $(0)^*_{\Coh{d}{R}}$ 
and a generic hypersurface 
intersection $\overline{R}$ of $R$. Here we use 
the characteristic $p$ version of Flenner's Bertini theorem.
Then finally we consider the non-vanishing of $\Coh{d-1}{R}$
in terms of non-vanishing of $(0)^*_{\Coh{d}{R}}$.

\section{Graded Isolated non $F$-rational Singularities}

This section summerizes the basic definitions and results concerning
isolated non $F$-rational singularities.
See, for example, \cite{HuSIX, HuAMS, KS8} for the detail.

\begin{Definition}
A sequence $x_1,\ldots, x_i$ in a ring $R$ is called 
{\em parameters} if for every $P\in\Spec{R}$ containing 
the sequnece, the image of the sequence in $R_P$ 
forms a part of sop. In other words, $\height((x_1,\ldots, x_i))=i$.
An ideal generated by a set of parameters is called 
a {\em parameter ideal}.
\end{Definition}

\begin{Definition}
For an ideal $I$ of a ring $R$, the {\em tight closure} $I^*$ 
of $I$ is defined by $I^* =\{r\in R\;\vert\; \mbox{there exists }c\in R^o
\mbox{ such that } cz^{q}\in I^{[q]}\mbox{ for all }q = p^e\gg 0\}$,
where $R^o$ is the set of elements outside the 
union of mininal primes of $R$ and 
$I^{[q]}= (r^q \;\vert\; r\in I)$.
\end{Definition}

\begin{Definition}
A Noetherian ring $R$ of prime characteristic is said to be {\em
$F$-rational} if every parameter ideal of $R$ is tightly closed.
\end{Definition}

We recall the basic properties of $F$-rational rings needed in
this paper. 

\begin{Proposition}[cf. \cite{HuAMS}]
\label{FratBasics}
Let $R$ be a Noetherian ring that 
is the homomorphic image of a Cohen-Macaulay ring. We have 
the following:
\begin{enumerate}
\item If $R$ is $F$-rational, then $R$ is Cohen-Macaulay.
\item A localization of an $F$-rational ring is $F$-rational.
\item If $R$ is regular, then $R$ is $F$-rational.
\end{enumerate}
\end{Proposition}

\begin{Definition}
Let $K$ be a field. Then a finitely generated $\NN$-graded $K$-algebra
$(R,\mm)$ is said to be a {\em generalized
Cohen-Macaulay} if one of the following equivalent conditions holds:
\begin{enumerate}
\item $\Proj(R)$ is an equidimensional Cohen-Macaulay projective scheme,
\item $R_P$ is an equidimensional 
Cohen-Macaulay local ring for every $P\in\Spec{R}- \{\mm\}$,
\item $\Coh{i}{R}$ is a finite length module for all $i<\dim R$.
      In particular, there exists $N_i\in\ZZ$ such that 
      $[\Coh{i}{R}]_n =0$ for all $n<N_i$.
\end{enumerate}
\end{Definition}

\begin{Definition}
Let $(R,\mm)$ be a finitely generated 
$\NN$-graded $K$-algebra or  a Noetherian local ring with
the maximal ideal $\mm$. Then we say that $R$ is an
{\em isolated non $F$-rational singularity} if 
$R_P$ is $F$-rational for all primes $P(\ne \mm)$.
\end{Definition}

In this paper, when we consider an isolated non $F$-rational
singularity, we always assume that it is equidimensional.

\begin{Definition}
An element $c\in R^o$ is called a {\em parameter test element}
if, for arbitrary ideal  $I$ generated by an sop and 
arbitrary $x\in I^*$, we have 
$c x^q \in I^{[q]}$ for all $q=p^e$.
\end{Definition}

In this definition, we can replace the ideal $I$ 
by any parameter ideal, including $I = (0)$, because 
of the following result, which is a special case 
of Exercise~2.12 \cite{HuAMS}.

\begin{Proposition}
Let $(R,\mm)$ be a local ring and assume that 
$c\in R^o$ is a test element for parameter ideals 
generated by sops. Then $c$ is a test element for 
arbitrary parameter ideal $I=(x_1,\ldots, x_i)$, $i<d=\dim R$,
including the case of $i=0$, i.e., $I=0$.
\end{Proposition}

Every isolated non $F$-rational singularity has 
an $\mm$-primary parameter test ideal. Namely,
\begin{Proposition}
\label{thm:VzCor}
Let $(R,\mm)$ be a reduced isolated non-$F$-rational singularity
of $\chara{R}=p>0$.
Then there exists an $\mm$-primary parameter test ideal $J\subset R$, i.e.,
every element $a\in J$ is a test element for any parameter ideal
in $R$.
\end{Proposition}
\begin{proof}
This follows immediately from  (3.9) \cite{Vz}.
See also Exer.~2.12 \cite{HuAMS}.
\end{proof}

By Prop.~\ref{FratBasics}, we know that an isolated 
non $F$-rational singularity  is generalized Cohen-Macaulay.
Also, in view of Prop.~\ref{thm:VzCor}, 
we will always consider reduced  rings.

\section{Lower local cohomologies of isolated non $F$-rational singularities}

This section gives a representation of lower local cohomologies of
isolated non $F$-rational singularities in terms of tight closure and
limit closure of unconditioned strong d-sequences (USD-sequences),
which is a refinement of the representation given by 
P.\ Schenzel and K.\ E.\ Smith.

\subsection{equidimensional hull, limit closure and germ closure}

In this subsection, we summerize some of the definitions and results
needed to give our representation of lower local cohomologies.

\begin{Definition}
A sequence of elements $x_1,\ldots, x_n$ in a commutative ring 
is said to be a {\em $d$-sequence} if for every $0\leq i\leq n-1$
and $k>i$ we have 
\begin{equation*}
   (x_1,\ldots, x_i) : x_{i+1}x_k = (x_1,\ldots, x_i):x_k.
\end{equation*}
A sequence $x_1,\ldots, x_n$ is said to be a {\em strong $d$-sequence}
if $x_1^{m_1},\ldots, x_n^{m_n}$ is a $d$-sequence for 
every arbitrary $m_i\geq 1$, $i=1,\ldots, n$. Finally, a sequence $x_1,\ldots, x_n$ is said to be
a {\em USD-sequence} (unconditioned strong d-sequence)
if every permutation of it is a strong 
$d$-sequence.
\end{Definition}

\begin{Definition}
Let $(R,\mm)$ be a Noetherian local ring with $d=\dim R$.
A system of parameters $x_1,\ldots, x_d$ is called 
{\em standard} if 
\begin{equation*}
(x_1,\ldots, x_d)\Coh{j}{R/(x_1,\ldots, x_{i-1})}=0
\qquad\mbox{for all }0\leq i+j\leq d,\; i\geq 1 
\end{equation*}
where we set $(x_1,\ldots, x_{i-1})=0$ if $i=1$.
\end{Definition}
Notice that for a standard sop we have
\begin{equation*}
   (x_1,\ldots, x_d) \Coh{i}{R}=0 \quad\mbox{for all }i< d.
\end{equation*}
We also note that the empty sequence is trivially a 
USD and standard sequence.
The notions of USD-sequence and standard sequence are actually
equivalent:
\begin{Proposition}[cf. (3.8)~\cite{Sch}]
\label{Sch3.8}
For an sop $x_1,\ldots, x_d$ in a Noetherian local ring $(R,\mm)$,
the following are equivalent:
\begin{enumerate}
\item [$(i)$] $x_1,\ldots, x_d$ is a standard system of parameters;
\item [$(ii)$] $x_1,\ldots, x_d$ is an unconditioned strong d-sequence
(USD-sequence).
\end{enumerate}
\end{Proposition}

The USD or standad property of the sequences
is preserved by hypersurface intersection. Namely,

\begin{Proposition}[cf. (3.2)~\cite{Sch}]
\label{hypersurfaceUSD}
Let $x_1,\ldots, x_i$ be a standard sequence in a Noetherian 
local ring $(R,\mm)$. Then the image of $x_1,\ldots, x_{i-1}$
in $\overline{R}:= R/x_iR$ is again a standard sequence of  
$\overline{R}$.
\end{Proposition}

For the existence of USD-sequence, we have 
\begin{Proposition}[(6.19) \cite{GY}]
\label{existenceUSD}
Let $(R,\mm)$ be a generalized Cohen-Macaulay local ring.
Then there exists $t\in\NN$ such that any system of parameters
in $\mm^t$ is  a USD-sequence.
\end{Proposition}
This implies that we can always obtain a USD-sequence by
taking hight enough power of an sop:
\begin{Corollary}
\label{usdExistence}
Let $(R,\mm)$ be a generalized Cohen-Macaulay local ring
and let $x_1,\ldots, x_d$, $d=\dim R$, be any system of 
parameters. Then $x_1^n, \ldots, x_d^n$ is a USD-sequence
for $n\gg 0$.
\end{Corollary}

\begin{Proposition}[cf. (5.11) \cite{HuSIX}]
\label{pteISusd}
Let $x_1,\ldots, x_i\in R$ be parameters which are also parameter 
test elements. Then $x_1,\ldots, x_i$ are a USD-sequence.
\end{Proposition}

\begin{Definition}
Let $R$ be a Noetherian ring.
For an ideal $I\subset R$, consider the minimal
primary decomposition $I = \Sect_iQ_i$. Then we define
$I^{umn} = \Sect Q_i$, where $Q_i$ runs over the primary components 
such that $\dim R/Q_i =  \dim R/\sqrt{Q_i} = \dim R/I$.
We call $I^{umn}$ the {\em unmixed hull} of $I$.
\end{Definition}

\begin{Definition}
For a parameter ideal $I=(x_1,\ldots, x_i)$ of a ring $R$,
we define the {\em limit closure} of $I$ as follows.
If $i\geq 0$, we define
\begin{eqnarray*}
   I^{\lim} 
    &=& \{z\in R \;\vert\;  {\bf x}^{s-1}z
    \in (x_1^s, x_2^s, \ldots, x_i^s) \mbox{ for some $s\in\NN$}\} \\
    &=&  \{z\in R \;\vert\;  {\bf x}^{s-1}z
    \in (x_1^s, x_2^s, \ldots, x_i^s) \mbox{ for all $s\gg 0$}\}\\  
    &=& \Union_{s=1}^\infty
       (x_1^s,\ldots, x_i^s): {\bf x}^{s-1}
\end{eqnarray*}
where ${\bf x} = {x_1}\cdots{x_i}$.
For $i=0$, we define $(0)^{\lim} = 0$.
\end{Definition}


\begin{Proposition}[cf. (2.5)~\cite{HuSm}]
\label{limclo1} Let $(R,\mm)$ be a Noetherian local ring and 
$I=(x_1,\ldots, x_i)$ $(i\geq 0)$ be a parameter ideal.
Then the following are equivalent:
\begin{enumerate}
\item [$(i)$] $z\in I^{\lim}$
\item [$(ii)$] an element 
$\eta = \displaystyle{\left[\frac{z}{x_1\cdots x_i}\right]}
\in \eCoh{I}{i}{R}$ is $0$.
\end{enumerate}
\end{Proposition}
\begin{proof}
Clear from the $\check{{\rm C}}$ech complex representation of 
local cohomology. 
\end{proof}

\begin{Proposition}[(5.4) \cite{HuSIX}]
Let $R$ be  equidimensional and the homomorphic image of 
a Cohen-Macaulay ring. If  $I\subset R$ is a parameter ideal,
then we have $I^{\lim}\subseteq I^*$. 
\end{Proposition}

\begin{Proposition}[(5.8) \cite{HuSIX}]
\label{prop5.8}
Let $(R,\mm)$ be  an equidimensional graded Noetherian  ring over 
a field with an $\mm$-primary parameter test ideal. 
Then for a parameter ideal $I\subset R$ 
such that $\height{I}<\dim R$, we have $I^{umn}=I^*$.
\end{Proposition}

\begin{Corollary}
\label{prop5.8Cor}
Let $(R,\mm)$ be an $\NN$-graded isolated non $F$-rational singularity
and let $I=(x_1,\ldots, x_i)$ be a parameter ideal such that 
$\height(I)< \dim R$. Then we have $I^{unm} = I^*$.
\end{Corollary}
\begin{proof}
Immediate from Prop.~\ref{thm:VzCor} and Prop.~\ref{prop5.8}.
\end{proof}

\begin{Definition}
Let $I = (x_1,\ldots, x_i)$ be a parameter ideal.
Then we define the 
{\em germ closure} of $I$ as follows. For $i\geq 1$, we define
\begin{equation*}
   I^{germ} = \sum_{j=1}^i (x_1,\ldots, \hat{x}_{j},,\ldots, x_i)^*.
\end{equation*}
Also, we define  $(0)^{germ} = 0$.
\end{Definition}

In the above definition, we note that $I^{germ} = (0)^*$ 
if $i=1$ and we have $I \subset I^{germ}$ for $i\ne 1$.
Also we always have $I^{germ}\subset I^*$.

The following proposition is a slight modification of Theorem~5.12
\cite{HuSIX}. The only difference is that, for $I = (x_1,\ldots, x_i)$
for $i\leq d$, only the case of $i=d$ is considered in \cite{HuSIX}. But
with almost the same proof we have  the following extension.

\begin{Proposition}[cf. (5.12) \cite{HuSIX}]
\label{husix:th5.12}
Let $(R,\mm)$ be an equidimensional local ring that is the 
homomorphic image of a Cohen-Macaulay local ring. 
Let $I = (x_1,\ldots, x_i)$, where $0\leq i\leq d$ and  
$x_1,\ldots, x_d$ is an sop which  are parameter test elements
Then we have $I^{germ} + I = I^{\lim}$,
in particular if $i\ne 1$ we have $I^{germ}=I^{\lim}$.
\end{Proposition}
\begin{proof}
We only show the case of $i=1$.  Other cases are 
similar to (5.12) \cite{HuSIX}.
Let $x\in R$ be a parameter test element. We have 
$(x)^{germ}= (0)^*_R$ by the definition. Let $z\in (x)^{\lim}$
be arbitrary. Then $x^sz \in (x^{s+1})$ for some $s\in\NN\cup \{0\}$,
so that $x^s(z - rx)=0$ for some $r\in R$. Assume that 
$z-rx\ne 0$. Since $x$ is a parameter test element, which 
also implies to be a test element for $(0)$, we have 
$z-rx \in (0)^*$ so that $z\in (x)+ (0)^* = (x)+ (x)^{germ}$.
Thus we have  $(x)^{\lim}\subset (x) + (x)^{germ}$.
For the converse inclusion, $(x)\subset (x)^{\lim}$ is immediate.
Again, since $x$ is a parameter test element, we have,
for $z\in (0)^*_R = (x)^{germ}$, $xz^q=0$ for any $q=p^e$.
In particular, $xz = 0\in (x^2)$, which implies $z\in (x)^{\lim}$.
Consequently, we have $(x)^{\lim}\supset (x) + (x)^{germ}$.
\end{proof}

\subsection{representation of lower local cohomologies}
\begin{Proposition}[cf. (6.8) and (6.9) \cite{KS8}]
\label{schenzelModpre}
Let $(R,\mm)$ be a Noetherian local ring that is
equidimensional, the homomorphic image of a Cohen-Macaulay local ring
and $R_P$ is Cohen-Macaulay for all primes $P\ne\mm$. Then there 
exists an sop $x_1,\ldots, x_d$, ($d=\dim R$) such that 
for $0\leq i<d$ we have
\begin{equation*}
   \Coh{i}{R} 
\iso
\frac{(x_1,\ldots, x_i)^{unm}}
{\sum_{j=1}^i(x_1,\ldots, \hat{x}_j,\ldots, x_i)^{unm}
 + (x_1,\ldots, x_i)}
\end{equation*}
Moreover if $i\ne 1$ we have
\begin{equation*}
   \Coh{i}{R} 
\iso
\frac{(x_1,\ldots, x_i)^{unm}}
{\sum_{j=1}^i(x_1,\ldots, \hat{x}_j,\ldots, x_i)^{unm}}
\end{equation*}
In fact, any sop $x_1,\ldots, x_d$ in $\mm^N$ for $N\gg 0$ 
has such a property.
\end{Proposition}

The following lemma, which we will use later,
is also used to show Prop.~\ref{schenzelModpre}.

\begin{Lemma}[6.9 \cite{KS8}]
\label{KS9-6.9}
In ths situation of Prop.~\ref{schenzelModpre},
we have $(x_1,\ldots, x_i): x_{i+1} = (x_1,\ldots, x_i)^{unm}$
for $0\leq i<d$, by replacing the sequence $x_1,\ldots, x_d$
by high powers, if necessary.
\end{Lemma}

\begin{proof}[Proof of Prop.~\ref{schenzelModpre}]
By Schenzel's formula (3.3)~\cite{Sch}, we have 
\begin{equation*}
\Coh{i}{R}\iso 
\frac{(x_1,\ldots, x_i):x_{i+1}}
     {\sum_{j=1}^i(x_1,\ldots, \hat{x}_j,\ldots, x_i): x_j
     +  (x_1,\ldots, x_i)}
\end{equation*}
for $i=0,\ldots, d-1$.
Thus we obtain the desired result by Lemma~\ref{KS9-6.9}.
Notice that for $i=2,\ldots, d-1$, we have
$(x_1,\ldots, x_i)
\subset \sum_{j=1}^{i}(x_1,\ldots, \hat{x}_j,\ldots, x_i):x_j$.
\end{proof}

Notice that the equidimensional hull of a graded ideal is again
graded. Then, by considering a graded version of Prop.~3.3 in
\cite{Sch}, which Prop.~\ref{schenzelModpre} bases on, we easily deduce the
following.

\begin{Corollary}
\label{schenzelMod}
Assume that $(R,\mm)$ is generalized Cohen-Macaulay. 
Then there exists a homogeneous sop $x_1,\ldots, x_d$, ($d=\dim R$) such
 that 
for $0\leq i < d$ we have 
\begin{equation*}
   \Coh{i}{R} 
\iso 
\frac{(x_1,\ldots, x_i)^{unm}}
{\sum_{j=1}^i(x_1,\ldots, \hat{x}_j,\ldots, x_i)^{unm} + (x_1,\ldots, x_i)}
(\delta_i)
\end{equation*}
where $\delta_i = \sum_{j=1}^{i}\deg(x_j)$.
Moreover if $i\ne 1$ we have
\begin{equation*}
   \Coh{i}{R} 
\iso
\frac{(x_1,\ldots, x_i)^{unm}}
{\sum_{j=1}^i(x_1,\ldots, \hat{x}_j,\ldots, x_i)^{unm}}(\delta_i).
\end{equation*}
\end{Corollary}

\begin{Theorem}
\label{1}
Let $(R,\mm)$ be an $\NN$-graded Noetherian $K$-algebra of dimension
$d=\dim R$ and $\chara{K}=p>0$.  Assume that $R$ is an equidimensional 
isolated non $F$-rational singularity. Then for any homogeneous
sop $x_1,\ldots, x_d \in\mm^N$ with $N\gg 0$, we have 
\begin{equation*}
\Coh{i}{R} \iso 
    \frac{(x_1,\ldots, x_i)^*}
         {(x_1,\ldots, x_i)^{\lim}}(\delta_i)
\iso 
\frac{(x_1^n,\ldots, x_i^n)^*}
     {(x_1^n,\ldots, x_i^n)^{\lim}}(n\delta_i)
\end{equation*}
for $i=0,\ldots, d-1$ and $n\in\NN$,
where we define $\delta_i =\sum_{j=1}^{i} \deg x_j$.
\end{Theorem}
\begin{proof}
Let $x_1,\ldots, x_d$ be a homogeneous sop. Since $R$ is  generalized
Cohen-Macaulay, Cor.~\ref{usdExistence} allows us to assume
that it is a USD-sequence by replacing it, if necessary,  by higher powers.
For $i=0,\ldots, d-1$, we compute
\begin{eqnarray*}
\Coh{i}{R} & \iso &
      \frac{(x_1,\ldots, x_i)^{unm}}
           {\sum_{j=1}^i(x_1,\ldots, \hat{x}_j,\ldots, x_i)^{unm}
              + (x_1,\ldots, x_i)}(\delta_i) 
	   \quad\mbox{by Cor.\ref{schenzelMod}} \\
           & \iso &
      \frac{(x_1,\ldots, x_i)^{*}}
           {(x_1,\ldots, x_i)^{germ} + (x_1,\ldots, x_i)}(\delta_i)
                    \quad\mbox{by Cor.\ref{prop5.8Cor}} \\
           & \iso &
      \frac{(x_1,\ldots, x_i)^{*}}
           {(x_1,\ldots, x_i)^{lim}}(\delta_i)
                 \quad\mbox{by Prop.~\ref{husix:th5.12}}\\
\end{eqnarray*}
The isomorphism
\begin{equation*}
\Coh{i}{R} \iso 
\frac{(x_1^n,\ldots, x_i^n)^*}
     {(x_1^n,\ldots, x_i^n)^{\lim}}(n\delta_i)
\qquad (n\in\NN)
\end{equation*}
is immediate since, for a USD-sequence $x_1,\ldots, x_i$,
its power $x_1^n,\ldots, x_i^n$ $(n\in\NN)$ is also a USD-sequence.
\end{proof}

As an immediate application of this representation, we recover 
some of the known results in a slightly general way.

\begin{Corollary}[cf. (2.7) \cite{HuSm}]
Let $(R,\mm)$ be an equidimensional isolated non $F$-rational singularity
and let $I=(x_1,\ldots, x_i)$, $i<d=\dim R$,
be a homogeneous parameter ideal such that 
with $I \subset \mm^\ell$ for $\ell\gg 0$.
Then for an integer $N_i$ the following are equivalent
\begin{enumerate}
\item [$(i)$] $[\Coh{i}{R}]_n=0$ for all $n<N_i$.
\item [$(ii)$] $I^* \subseteq I^{\lim} + R_{\geq \delta_i + N_i}$ 
where $\delta_i = \sum_{j=1}^{i}\deg(x_j)$.
\end{enumerate}
\end{Corollary}

In a special case, we can give a lower bound of the vanishing degree
$N_i$. For an $\NN$-graded module $M$, we set $b(M) = \min\{i \;\vert\;
M_i\ne 0 \}$. 
%
%
If $(R, \mm)$ is a generalized Cohen-Macaulay local ring, 
we have 
\begin{equation*}
[\Coh{i}{R}]_n=0\qquad\mbox{for}\quad n< -it,\quad i<\dim R,
\end{equation*}
if there exists a homogeneous sop $x_1,\ldots,x_d$ 
as in Th.~\ref{1} with $t = \min\{\deg x_1,\ldots, \deg x_d\}$
((2.4)\cite{HM95}).
The following result
gives a refinement for isolated non $F$-rational singularities.

\begin{Corollary}
\label{vanishingDeg}
Let $K$ be an algebraically closed field of $\chara{K}=p>0$ and 
let $(R,\mm)$ be a Noetherian standard $\NN$-graded domain  with $R_0=K$ 
and $\mm = \Dirsum_{n>0}R_n$.  Assume that $R$ is an isolated 
non $F$-rational singularity and let 
\begin{equation*}
t = \min\{ N\in\NN \;\vert\; {\bf x}= x_1,\ldots, x_d\in\mm^N 
        \mbox{ where ${\bf x}$ is as in Th.~\ref{1}}\}.
\end{equation*}
Then for all $i<d=\dim{R}$ and for all $n< -(i-1)t$ we have 
$[\Coh{i}{R}]_n=0$.
\end{Corollary}
\begin{proof}
By Th.~\ref{1} we can find an sop $x_1,\ldots, x_d\in\mm^N$ such 
that for every $i<d$
\begin{equation*}
   \Coh{i}{R}\iso \frac{(x_1,\ldots,x _i)^*}{(x_1,\ldots,
    x_i)^{\lim}}(\delta_i)
\qquad \mbox{where}\quad
\delta_i = \sum_{j=1}^{i}\deg(x_j)
\end{equation*}
The lower bound of the degrees of the non-zero elements in 
$(x_1,\ldots, x_i)$ is\\
 $\min\{\deg x_1,\ldots, \deg x_i\}$, 
so that we have 
$b((x_1,\ldots, x_i)^*(\delta_i)) 
= b((x_1,\ldots, x_i)^*) - \delta_i$,
which is, by Exer.~4.5 and ~4.6 in \cite{HuAMS},
equal to 
$b((x_1,\ldots, x_i)) - \delta_i$\\
$= \min\{\deg{x_1},\ldots, \deg{x_i}\} - \delta_i
\geq -(i-1)t$ as required.
\end{proof}

For another application, 
we can consider,
under the condition of Theorem~\ref{1},
the natural map $\rho_i : \Coh{i}{R} \To \Coh{i+1}{R}$
for $i=0,\ldots, d-2$ defined by $\rho_i(r + (x_1,\ldots, x_i)^{\lim})
= r + (x_1,\ldots, x_{i+1})^{\lim}$ with $r\in (x_1,\ldots, x_i)^*$.
For the degree preserving homomorphism, we 
can also consider 
$\psi_i(r + (x_1,\ldots, x_i)^{\lim})
= rx_{i+1} + (x_1,\ldots, x_{i+1})^{\lim}$.
Unfortunately, they are all trivial.

\begin{Corollary}
\label{0map}
The maps $\rho_i$ and $\psi_i$ are  trivial.
\end{Corollary}
\begin{proof}
We first show that 
$(x_1,\ldots, x_i):x_{i+1} \subset (x_1,\ldots, x_{i+1})^{\lim}$
for $i<d$.  Take any  $r + (x_1,\ldots, x_i)^{\lim}\in
\displaystyle{\frac{(x_1,\ldots, x_i)^*}{(x_1,\ldots, x_i)^{\lim}}}
\iso\Coh{i}{R}$.
By Prop.~\ref{prop5.8} and Lemma~\ref{KS9-6.9}, we have 
$r\in (x_1,\ldots, x_i)^* = (x_1,\ldots, x_i):x_{i+1}$,
i.e., $x_{i+1}r \in (x_1,\ldots, x_i)$,
so that, for any $s\in\NN$, we have 
$x_1^s\cdots x_i^s x_{i+1}^sr 
\in (x_1^{s+1},\ldots, x_{i+1}^{s+1})$,
which implies that $r \in (x_1,\ldots, x_{i+1})^{\lim}$.
Hence
$\rho_i=0$.
$\psi_i=0$ is similarly proved.
\end{proof}

\section{Tight Closure of Zeros and Isolated Singularities}

In this section, we consider the question of how the tight
closure of zero in the highest local cohomology 
$\Coh{d}{R}$ controls the vanishing of the lower local
cohomologies of isolated singularities.

\begin{Definition}[tight closure of zero]
For a $R$-module $M$, we define 
\begin{equation*}
(0)^*_{M}
= \{x \in M \;\vert\; \mbox{there exists $c\in R^{o}$ such that 
$c x^q =0$ in $M$ for $q=p^e\gg 0$}  \}
\end{equation*}
and call it {\em tight closure of $0$} in $M$.
\end{Definition}
We are particularly interested in the tight closure of 
zero in the highest local cohomology $\Coh{d}{R}$, $d=\dim R$.
We give a graded version of the result by Goto and Nakamura.

\begin{Proposition}
\label{tc0}
Let $(R,\mm)$ be an $\NN$-graded isolated non $F$-rational singularity.
Suppose that a homogeneous
sop ${\bf x} = x_1,\ldots, x_d$ forms a USD-sequence. Then,
by replacing ${\bf x}$ by a higher power ${\bf x}^N$ $(N>0)$ if necessary,
we have 
\begin{equation*}
   (0)^*_{\Coh{d}{R}}
    = \Union_{n>0}Z_n(n\delta)
\qquad\mbox{where}\quad
   Z_n = \frac{(x_1^n,\ldots, x_d^n)^*}{(x_1^n,\ldots, x_d^n)^{\lim}},
\end{equation*}
where $\delta = \sum_{j=1}^d\deg(x_j)$.
\end{Proposition}
\begin{proof}
By considering the graded version of Prop.~2.1 \cite{GN02}, we have 
\begin{equation*}
(0)^*_{\Coh{d}{R}}
    = \Union_{n>0}Z'_n(n\delta)
\quad\mbox{where}\quad
   Z'_n = \frac{(x_1^n,\ldots, x_d^n)^*}
          {\sum_{j=1}^n(x_1^n,\ldots, \hat{x}^n_j\ldots, x_d^n): x_j^n
	        + (x_1^n, \ldots, x_d^n)}.
\end{equation*}
Then by Cor.~\ref{prop5.8Cor}, 
Prop.~\ref{husix:th5.12} and 
Lemma~\ref{KS9-6.9}, we obtain the desired result.
\end{proof}

Notice that, for $i<d=\dim R$, we have a natural map 
\begin{equation*}
\Coh{i}{R}
\iso \frac{(x_1,\ldots, x_i)^*}{(x_1,\ldots, x_i)^{\lim}}(\delta_i)
\To \frac{(x_1,\ldots, x_d)^*}{(x_1,\ldots, x_d)^{\lim}}(\delta_d)
\subset (0)^* \subset \Coh{d}{R}
\end{equation*}
by sending $r + (x_1,\ldots, x_i)^{\lim}$ to 
 $r + (x_1,\ldots, x_d)^{\lim}$
or as a degree preserving map
$rx_{i+1}\cdots x_{d} + (x_1,\ldots, x_d)^{\lim}$. 
But they are actually  $0$-maps, which can be
proved similarly to Cor.~\ref{0map}.

Now we consider relation between the tight closures of zero
of $\Coh{d}{R}$ and $\Coh{d-1}{\overline{R}}$, where 
$\overline{R}$ is a hypersurface intersetion of $R$.
It is well known that Flenner's Bertini theorem  also 
holds for positive characteristics, which is, as far as 
the author is concerned, not stated in the literature. 
We give the statement for the readers convenience.
\begin{Proposition}[cf. Satz~(4.1) \cite{Fl}]
\label{charPsatz4.1}
Let $k$ be an infinite field and $(R,\mm)$ be a Noetherian local
$k$-algebra, whose residue class ring $K:=R/\mm$ is separable 
over $k$. Let $x_1,\ldots, x_d\in\mm$ be such that
$I = (x_1,\ldots, x_d)$ is an $\mm$-primary ideal.
Assume that $R_P$ is regular for all $P\in \Spec{R}-\{\mm\}$.
Then for a general linear form
$x_\alpha = \sum_{i=1}^d \alpha_i x_i$ with $\alpha \in k^d$
and any $P\in \Spec{R/x_\alpha}- \{\mm'\}$ where $\mm' = \mm/x_\alpha R$, 
$(R/x_\alpha R)_P$ is again regular.
\end{Proposition}

Let $(R,\mm)$ be an isolated singularity with $d=\dim R$. 
Then, by Prop.~\ref{hypersurfaceUSD}, ~\ref{charPsatz4.1} and~\ref{tc0}
together with Cor.~\ref{usdExistence},
we can choose a homogeneous sop $x_1,\ldots, x_d$, which is a 
USD-sequence,  such that 
$\overline{R} = R/x_dR$ is also an isolated singularity 
and 
\begin{equation*}
(0)^*_{\Coh{d}{R}} 
 = \Union_{n>0}\displaystyle{
 \frac{(x_1^n,\ldots, x_d^n)^*}{(x_1^n,\ldots, x_d^n)^{\lim}}(n\delta_d)
 }
\quad\mbox{and}\quad
(0)^*_{\Coh{d-1}{\overline{R}}} 
 = \Union_{n>0}\displaystyle{
 \frac{(\overline{x}_1^n,\ldots, \overline{x}_{d-1}^n)^*}
{(\overline{x}_1^n,\ldots, \overline{x}_{d-1}^n)^{\lim}}(n\delta_{d-1})
 }
\end{equation*}
where $\delta_j = \sum_{k=1}^j\deg(x_k)$
and $\overline{x}_j$ is the image of $x_j$ in $\overline{R}$.
Now we can prove
\begin{Proposition}
\label{mainSub}
Let $(R,\mm)$ be an $\NN$-graded isolated non $F$-rational 
singularity over the field $K$ of $\chara{K}=p>0$ with $R_0=K$. 
Let $d=\dim R$.
Assume that $x_1,\ldots, x_d$ is a homogeneous sop
and consider the natural surjection
\begin{equation*}
  \varphi : R \To \overline{R} = R/x_dR,
\end{equation*}
whose image $\varphi(r)$ will be denoted by $\overline{r}$.
Then,
\begin{enumerate}
\item [$(i)$]
For any $n\in\ZZ$ we have 
\begin{equation*}
\mbox{if}\quad \left[
   \frac{(x_1,\ldots, x_{d})^*}
    {(x_1,\ldots, x_{d})^{\lim}}
  \right]_n \ne 0,
\quad
\mbox{then}\quad
  \left[
  \frac{(\overline{x}_1,\ldots, \overline{x}_{d-1})^*}
   {(\overline{x}_1,\ldots, \overline{x}_{d-1})^{\lim}}
  \right]_n\ne 0,
\end{equation*}
\item [$(ii)$] Assume also that 
$x_1,\ldots, x_d$ is a USD-sequence and that 
$R/x_d^\ell R$ is a non $F$-rational singularity
for all $\ell\in\NN$. 
If $[(0)^*_{\Coh{d}{R}}]_n\ne 0$ for $n\in\ZZ$, then 
there exists $\ell\in\NN$ such that 
$[(0)^*_{\Coh{d-1}{R/x_d^\ell R}}]_{n+a}\ne 0$,
where $a = \ell\cdot \deg(x_d)$.
\end{enumerate}
\end{Proposition}
\begin{proof}
We first show $(i)$.
By persistence of tight closure (cf. Th.~2.3 \cite{HuAMS}) we have 
$\varphi((x_1,\ldots, x_d)^*) \subset (\overline{x}_1,\ldots, \overline{x}_d)^*= (\overline{x}_1,\ldots, \overline{x}_{d-1})^*$.
Thus 
$\varphi$ induces the homomorphism
\begin{equation*}
 \overline{\varphi}: 
   (x_1,\ldots, x_d)^*
 \To
   \frac{(\overline{x}_1,\ldots, \overline{x}_{d-1})^*}
   {(\overline{x}_1,\ldots, \overline{x}_{d-1})^{\lim}}.
\end{equation*}
We now show $\Ker\overline{\varphi}\subset (x_1,\ldots, x_d)^{\lim}$.
Let $r\in\Ker\overline{\varphi}$ be arbitrary. Then
$\overline{r}:= 
\varphi(r)\in (\overline{x}_1,\ldots, \overline{x}_{d-1})^{\lim}$
and  we have 
$\overline{x}_1^s\cdots\overline{x}_{d-1}^s\cdot\overline{r}
\in (\overline{x}_1^{s+1}, \ldots, \overline{x}_{d-1}^{s+1})$
for some $s\in\NN\cup \{0\}$, which implies that $x_1^s\cdots x_{d-1}^s r
\in (x_1^{s+1},\ldots, x_{d-1}^{s+1}, x_d)$ and then
$x_1^s\cdots x_d^s r \in (x_1^{s+1},\ldots, x_d^{s+1})$. 
Hence $r\in (x_1,\ldots, x_d)^{\lim}$ as required. 
Consequently, we have the following diagram:
\begin{equation}
\label{maindiagram}
\frac{(x_1,\ldots, x_d)^*}{(x_1,\ldots, x_d)^{\lim}}
\overset{nat}{\longleftarrow}
\frac{(x_1,\ldots, x_d)^*}{\Ker\overline{\varphi}}
\overset{\hat{\varphi}}{\hookrightarrow}
\frac{(\overline{x}_1,\ldots, \overline{x}_{d-1})^*}
     {(\overline{x}_1,\ldots, \overline{x}_{d-1})^{\lim}}
\end{equation}
where $nat$ is the natural surjection and $\hat{\varphi}$ 
is the induced embedding. Taking the degree $n$ fraction of 
this diagram, we immediately know that $(i)$ holds.

Now we show $(ii)$. Assume that $[(0)^*_{\Coh{d}{R}}]\ne 0$.
Then by Prop.~\ref{tc0}  there exists $\ell\in\NN$ such that 
\begin{equation*}
0\ne 
\left[
   \displaystyle{
    \frac{(x_1^\ell,\ldots, x_d^\ell)^*}
         {(x_1^\ell,\ldots, x_d^\ell)^{\lim}}(\ell\cdot\delta_d)
   }
\right]_n 
=
\left[
   \displaystyle{
    \frac{(x_1^\ell,\ldots, x_d^ell)^*}
         {(x_1^\ell,\ldots, x_d^\ell)^{\lim}}
   }
\right]_{n+\ell\cdot \delta_d}.
\end{equation*}
Thus by $(i)$ with $x_1,\ldots, x_d$ replaced by
$x_1^\ell,\ldots, x_d^\ell$, we have 
\begin{equation*}
0\ne 
\left[
   \displaystyle{
    \frac{(\overline{x}_1^\ell,\ldots, \overline{x}_{d-1}^\ell)^*}
         {(\overline{x}_1^\ell,\ldots, \overline{x}_{d-1}^\ell)^{\lim}}
   }
\right]_{n+\ell\cdot \delta_d} 
=
\left[
   \displaystyle{
    \frac{(\overline{x}_1^\ell,\ldots, \overline{x}_{d-1}^\ell)^*}
         {(\overline{x}_1^\ell,\ldots, \overline{x}_{d-1}^\ell)^{\lim}}
          (\ell\cdot \delta_{d-1})
   }
\right]_{n+\ell\cdot \deg (x_d)} 
\end{equation*}
Then we obtain the required result by Prop.~\ref{tc0}. 
\end{proof}

Prop.~\ref{mainSub} assumes that hypersurface intersections 
of isolated non $F$-rational singularities are again isolated 
non $F$-rational singularities. However, $F$-rationality behaves 
rather badly with hypersurface intersections, since 
an $F$-rational local ring has negative $a$-invariant 
(see \cite{FW}) but hypersurface intersection increases $a$-invariant.
But for isolated singularity, which is a special case of 
isolated non $F$-rational singularity, hypersurface intersections 
behave well thanks to Bertini type theorems.

Now we come to the main result of this section.

\begin{Theorem}
\label{main}
Let $(R,\mm)$ be a reduced $\NN$-graded equidimensional isolated singularity
over the field $K$ of $\chara{K}=p>0$ with $R_0=K$, $d=\dim R$ and 
 $\sharp{K}=\infty$. 
Assume that $[(0)^*_{\Coh{d}{R}}]_n\ne 0$
for some $n\in\ZZ$ and consider 
a homogeneous regular element $x$, $a = \deg(x)$,
that is also a part of a USD-sequence.
Then 
\begin{enumerate}
\item [$(i)$] If 
the multiplication $[\Coh{d}{R}]_{n}\overset{x}{\To} [\Coh{d}{R}]_{n+a}$
is injective, then we have\\
 $[\Coh{d-1}{R}]_{n+\ell\cdot a}\ne 0$ 
for some $\ell\in\NN$.
\item Otherwise, if there is a USD sequence $x_1,\ldots, x_d$ with $x=x_d$
such that the degree $n+\delta_d$ fragment of  $(x_1,\ldots, x_d)^{\lim} -
\varphi^{-1}((\overline{x}_1,\ldots, \overline{x}_{d-1})^{\lim})$
is nonempty, where $\varphi : R \To R/x_dR$ is the natural surjection whose
image of $r\in R$ is denoted by $\overline{r}$.
Then we have  $[\Coh{d-1}{R}]_{n+a}\ne 0$.
\end{enumerate}
\end{Theorem}
\begin{proof}
Let $x_1,\ldots, x_d\in\mm$ be a USD-sequence,
whose existence is assured by Prop.~\ref{existenceUSD}.
Since $R$ is reduced, we have $\depth R >0$ and thus we 
can take $x_d$ general enough to be $R$-regular. 
Consider the 
short exact sequence:
\begin{equation*}
   0 \To R(-a) \overset{x_d}{\To} R \overset{\varphi}{\To}
    \overline{R}\To 0,
\qquad a:= \deg(x_d)
\end{equation*}
where we set $\overline{R}:= R/(x_d)R$ and 
$\varphi: R \To R/(x_d)R$ is the natural map. 
Also by Prop.~\ref{charPsatz4.1}, we can assume that 
$\overline{R}$ is again an isolated singularity.

Since $x_d$ is a part of the 
USD-sequence, it is standard by Prop.~\ref{Sch3.8},
so that we have  $x_d \Coh{d-1}{R}=0$. 
Thus, from the above short exact sequence,
we obtain the long exact sequence 
\begin{equation*}
0\To \Coh{d-1}{R}\To\Coh{d-1}{\overline{R}}
\overset{\psi}{\To} \Coh{d}{R}(-a)\overset{x_d}{\To} \Coh{d}{R}\To 0,
\end{equation*}
and we consider the degree $n+a$ fragment:
\begin{equation*}
0\To [\Coh{d-1}{R}]_{n+a} \To
     [\Coh{d-1}{\overline{R}}]_{n+a} \overset{\psi}{\To }
     [\Coh{d}{R}]_n \overset{x_d}{\To} 
     [\Coh{d}{R}]_{n+a} \To 0.
\end{equation*}

Now we show $(i)$. Assume that the multiplcation by $x_d$ is injective, 
we have an isomorphism $[\Coh{d-1}{R}]_{n+a} \iso 
[\Coh{d-1}{\overline{R}}]_{n+a}$.
By Prop.~\ref{tc0} we have the inclusion
\begin{equation*}
   \frac{(\overline{x}_1,\ldots, \overline{x}_{d-1})^*}
    {(\overline{x}_1,\ldots, \overline{x}_{d-1})^{\lim}}(\delta_{d-1})
    \subset
    (0)^*_{\Coh{d-1}{\overline{R}}}
    \subset \Coh{d-1}{\overline{R}}
\end{equation*}
so that, if 
\begin{equation*}
   \left[
   \frac{(\overline{x}_1,\ldots, \overline{x}_{d-1})^*}
    {(\overline{x}_1,\ldots, \overline{x}_{d-1})^{\lim}}(\delta_{d-1})
   \right]_{n+a} \ne 0,
\end{equation*}
then we have 
$[\Coh{d-1}{\overline{R}}]_{n+a} \iso [\Coh{d-1}{R}]_{n+a}\ne 0$.
Then $(i)$ follows  by Prop.~\ref{mainSub}(ii) by replacing,
if necessary, $x_d$ by $x_d^\ell$ for some $\ell\in\NN$.

Next we show $(ii)$. If $x_d\in \mm^N$ for $N\gg 0$, $x_d$ is 
a parameter test element by Prop.~\ref{thm:VzCor}.
By Prop~4.4 \cite{KS1}
we know $J \subset \Ann_R((0)^*_{\Coh{d}{R}})$
where $J$ is the parameter test ideal.
Thus $x_d : [\Coh{d}{R}]_n \overset{x_d}{\To} 
     [\Coh{d}{R}]_{n+a}$
is not injective if $[(0)^*_{\Coh{d}{R}}]_n\ne 0$.
Thus we have the exact sequence
\begin{equation*}
   0\To [\Coh{d-1}{R}]_{n+a} \To 
        [\Coh{d-1}{R}]_{n+a} 
	\overset{\psi}{\To}
	\Im\psi \To 0
\end{equation*}
where $0\ne [(0)^*_{\Coh{d}{R}}]_n \subset \Im\psi$. By 
considering the $\check{{\rm C}}$ech complex representation
of local cohomology, we know that,
for 
\begin{equation*}
\eta = \left[
	\frac{\overline{r}}{\overline{x}_1\cdots\overline{x}_{d-1}}
      \right]
\in \Coh{d-1}{\overline{R}}
\end{equation*}
we have 
\begin{equation*}
   \psi(\eta)  = \left[
		  \frac{r}{x_1\cdots x_d}
                \right]
   \in \Coh{d}{R}.
\end{equation*}
Also we know that 
\begin{eqnarray*}
\Ker\psi &=& \left\{
               \eta = \left[
		       \frac{\overline{r}}
		            {\overline{x}_1\cdots\overline{x}_{d-1}}
                     \right]
	       \;\vert\;
	       \left[
		\frac{r}{x_1\cdots x_d}
	       \right]=0\mbox{ in }\Coh{d}{R}
            \right\}\\
         & = & \left\{
               \eta = \left[
		       \frac{\overline{r}}
		            {\overline{x}_1\cdots\overline{x}_{d-1}}
                     \right]
	       \;\vert\;
	          r\in (x_1,\ldots, x_d)^{\lim}
            \right\}\\
\end{eqnarray*}
by Prop.~\ref{limclo1}.
We want to show that $\Ker\psi\ne 0$.
Now consider the diagram $(\ref{maindiagram})$ in the proof of 
Prop.~\ref{mainSub}. We consider the submodule
\begin{equation*}
L :=\Im\hat{\varphi} \subset 
\frac{(\overline{x}_1,\ldots, \overline{x}_{d-1})^*}
     {(\overline{x}_1,\ldots, \overline{x}_{d-1})^{\lim}}
\subset (0)^*_{\Coh{d-1}{\overline{R}}}(-\delta_{d-1}).
\end{equation*}
Then $\psi([L]_{n+\delta_d}) \subset [(0)^*_{\Coh{d}{R}}]_n$,
since 
\begin{equation*}
   \psi(\eta)  = \left[
		  \frac{r}{x_1\cdots x_d}
                \right]
   \in (0)^*_{\Coh{d}{R}}
\end{equation*}
if and only if $r\in(x_1,\ldots, x_d)^*$
(see for example \cite{HuSIX}). 
Now, if $r \in (x_1,\ldots, x_d)^{\lim} - \Ker\overline{\varphi}$,
then 
the image of $r$ in $(x_1,\ldots, x_d)^*/\Ker\overline{\varphi}$
in the diagram $(\ref{maindiagram})$
give the non-zero element 
$\eta \in (\overline{x}_1,\ldots, \overline{x}_{d-1})^*/
(\overline{x}_1,\ldots, \overline{x}_{d-1})^{\lim}$,
which contributes to a non-zero element in $\Coh{d-1}{R}$.
Moreover we have $\psi(\eta)=0$, 
since $\psi(\eta)$ is the image of $r$ in 
$(x_1,\ldots, x_d)^*/(x_1,\ldots, x_d)^{\lim}$. Thus 
$\Ker\psi\ne 0$ so that we must have $[\Coh{d-1}{R}]_{n+a}\ne 0$
as required.
\end{proof}

\begin{Remark}
{\em 
As we saw in the proof, the condition $(i)$ of Th.~\ref{main} is
satisfied only when $\deg(x)$ is small large enough.  On the other hand,
Prop.~\ref{existenceUSD} tells that, generally speaking, in order to
obtain a part of USD-sequence we need to choose elements of high degree.
Thus $(i)$ could be satisfied in relatively few cases.

Also, since $[H^{d-1}{R}]_k=0$ for $k\gg 0$,   Th.~\ref{main} 
means that $(x_1,\ldots, x_d)^{\lim}$ is 
the preimage of $(\overline{x}_1,\ldots, \overline{x}_{d-1})^{\lim}$
via the natural map $\varphi$ when $\deg(x_d)$ is large enough.
}
\end{Remark}

The following result describes the injectivity of the morphism
$nat$ in the diagram (\ref{maindiagram}).

\begin{Corollary}
\label{corollary}
Let $(R,\mm)$ be a reduced $\NN$-graded equidimensional isolated singularity
over the field $K$ of $\chara{K}=p>0$ with $R_0=K$,
$d=\dim R$ and $\sharp{K}=\infty$.
For an sop $x_1,\ldots, x_d \in \mm^N$ with $N\gg 0$
and for every $n\in\ZZ$ with
$[(0)^*_{\Coh{d}{R}}]_n\ne 0$, we have 
\begin{equation*}
\left[
 \displaystyle{
 \frac{(x_1,\ldots, x_d)^*}{(x_1,\ldots, x_d)^{\lim}}
 }
\right]_{n+\delta_d}
\subset 
\left[
 \displaystyle{
 \frac{(\overline{x}_1,\ldots, \overline{x}_{d-1})^*}
      {(\overline{x}_1,\ldots, \overline{x}_{d-1})^{\lim}}
 }
\right]_{n+\delta_d}
\end{equation*}
where $\overline{x}_j$ is the image of $x_j$ 
in $\overline{R}:= R/x_d R$.
\end{Corollary}
\begin{proof}
Since $[\Coh{d-1}{R}]_{n+a}=0$ for $a\gg 0$, we know that 
the conditions $(i)$ and $(ii)$ in Th.~\ref{main} do not hold 
for sop's with high enough degrees, which implies 
that $[(x_1,\ldots, x_d)^{\lim}]_{n+\delta_d} = 
[\varphi^{-1}((\overline{x}_1,
\ldots, \overline{x}_{d-1})^{\lim})]_{n+\delta_d}$ as required.
\end{proof}


\end{document}